# NOTE ON NEAR SUBNORMAL WEIGHTED SHIFTS[*]

WANG GONGBAO

(Department of Basic Courses, Naval University of Engineering, Wuhan 430033, China)

MA JIPU

( Department of Mathematics, Nanjing University, Nanjing 210093, China )

**ABSTRACT.** In this note, necessary and sufficient conditions are obtained for unilateral weighted shifts to be near subnormal . As an application of the main results, many answers to the Hilbert space problem 160 are presented at the end of the paper.

**MR(2000) Subject Classification.**  47B20,  47B37.

## 1.  INTRODUCTION

In [1], two operator classes of *D*-near subnormal operators and near subnormal operators were introduced. The elementary properties of such operators were obtained by use of the M-P generalized inverse. Moreover, a new necessary and sufficient condition for an operator to be subnormal was presented.

For convenience, some relevant concepts and results of [1] are cited as follows.

**Definition 1.** Let $A, D \in B(H)$ and $D \geq 0$. If there exists a constant $m > 0$ such that $D \geq m A^* D A$, then A is called a *D*-near subnormal operator.

**Definition 2.** Let *A* be a hyponormal operator and $Q_A = A^* A - A A^*$ ($Q_A \geq 0$). If A is a $Q_A$ − near subnormal operator, then *A* is said to be a near subnormal operator.

**Theorem.** Suppose *A* is a hyponormal operator, then *A* is near subnormal   if and only if

---

[*] Research supported by the National Natural Science Foundation of China (19971039) and the Science Foundation of Naval University of Engineering (2002706).

(1) $N(Q_A) \in LatA$ and (2) $Q_A^{\frac{1}{2}} A Q_A^{+\frac{1}{2}} \in B(H)$.

In the above theorem, $N(Q_A), LatA$ and $Q_A^{+\frac{1}{2}}$ denote the null space of $Q_A$, the invariant subspace lattice of $A$ and the M-P generalized inverse of $Q_A^{\frac{1}{2}}$ respectively(cf.[2]).

From [1], we know that the class of subnormal operators is a proper subset of near subnormal operators. And the class of near subnormal operators is a proper subset of hyponormal operators. In [3], Stampfli presented a necessary and sufficient condition for weighted shifts to be subnormal. In [4], we gave the same criterion with the result of [3] by use of the M-P generalized inverse. But the method was completely different from Stampfli's work. As an application of the theorem of [4], we also got many answers to the Hilbert space problem 160(cf. [5]).

The main purpose of this note aims at discovering which unilateral weighted shifts are near subnormal. Necessary and sufficient conditions are given in terms of the weight sequence. In the 160th problem of [5], it is required to find an operator which is hyponormal but not subnormal. Halmos regarded this problem as non-trivial. And there were a few answers to the question (cf. [1],[5],[6]and[7]). As an application of the main results of this note, we present a great many simple answers to it. Therefore, we solve the Hilbert space problem 160 further.

## 2. MAIN RESULTS

Let $H$ be a separable complex Hilbert space with an orthonormal basis $\{e_n\}_{n=1}^{\infty}$, and $B(H)$ the algebra of all bounded linear operators on $H$. Let $S$ be a unilateral weighted shift with weight sequence $\{a_n\}_{n=1}^{\infty}$, that is, $Se_n = a_n e_{n+1}$ for all $n \geq 1$, where $a_n \in \mathbf{C}$ and $\sup_{n \geq 1}|a_n| < +\infty$. Since $a_n = e^{i\theta_n}|a_n|$, $S$ is unitarily equivalent to the unilateral weighted shift with weight sequence $\{|a_n|\}_{n=1}^{\infty}$. Without real loss of generality, we may assume that every element of $\{a_n\}_{n=1}^{\infty}$ is non-negative.

It is well-known that $S$ is hyponormal if and only if its weight sequence $\{a_n\}_{n=1}^{\infty}$ is increasing. For a hyponormal unilateral weighted shift $S$, we may assume $a_n \neq 0$ for all

n≥1. In fact, if $a_1 = a_2 = \cdots = a_m = 0$, then $M = span\{e_1, e_2, \cdots, e_m\}$ reduces $S$ and $S$ is normal on $M$. For this reason, we assume all weight sequences to have no zeros throughout this paper.

Now, the main results of the note as follows.

**Theorem 1.** Let $S$ be a hyponormal operator with $Se_n = a_n e_{n+1}$ for all $n \geq 1$.

(1) If $0 < a_1 < a_2 < a_3 < \cdots$, then $S$ is near subnormal if and only if

$$\sup_{n \geq 1} \{a_n (\frac{a_{n+1}^2 - a_n^2}{a_n^2 - a_{n-1}^2})^{\frac{1}{2}}\} < +\infty, \text{ where } a_0 = 0.$$

(2) If $0 < a_1 \leq a_2 \leq a_3 \leq \cdots$, then $S$ is near subnormal if and only if there exists some integer $i \geq 1$ such that $0 < a_1 < a_2 < \cdots < a_{i-1} < a_i$ and $a_n = a_i$ for $n = i+1, i+2, \cdots$.

Proof. (1) Since $Se_n = a_n e_{n+1}$ for all $n \geq 1$, we have

$$Q_S = S^*S - SS^* = diag\{a_1^2, a_2^2 - a_1^2, a_3^2 - a_2^2, \cdots\}. \tag{*}$$

Therefore $Q_S^{\frac{1}{2}} = diag\{a_1, (a_2^2 - a_1^2)^{\frac{1}{2}}, (a_3^2 - a_2^2)^{\frac{1}{2}}, \cdots\}$ and the M-P generalized inverse of $Q_S^{\frac{1}{2}}$ is $Q_S^{+\frac{1}{2}} = diag\{\frac{1}{a_1}, \frac{1}{(a_2^2 - a_1^2)^{\frac{1}{2}}}, \frac{1}{(a_3^2 - a_2^2)^{\frac{1}{2}}}, \cdots\}$. By some calculations, it follows that $Q_S^{\frac{1}{2}} S Q_S^{+\frac{1}{2}}$ is also a unilateral weighted shift with weight sequence $\{a_n (\frac{a_{n+1}^2 - a_n^2}{a_n^2 - a_{n-1}^2})^{\frac{1}{2}}\}$.

From $(*)$, it is evident that $N(Q_S) = \{0\}$. Hence (1) holds by the theorem of Section 1.

(2) For the necessity, let $i$ be the smallest positive integer such that $a_i = a_{i+1}$. Thus $i \geq 1$ and $0 < a_1 < a_2 < \cdots < a_{i-1} < a_i$, $a_i = a_{i+1}$. On the other hand, it is easy to see that $e_{i+1} \in N(Q_S)$. Since $S$ is near subnormal, certainly $N(Q_S) \in LatS$. It follows that $Se_{i+1} = a_{i+1} e_{i+2} \in N(Q_S)$. So $e_{i+2} \in N(Q_S)$ for $a_{i+1} \neq 0$. Similarly, we have $\{e_n\}_{n=i+1}^{\infty} \subset N(Q_S)$. This implies that $a_n = a_i$ for $n = i+1, i+2, \cdots$.

Conversely, by the assumption it follows that $Q_S = diag\{a_1^2, a_2^2 - a_1^2, \cdots, a_i^2 - a_{i-1}^2, 0, 0, \cdots\}$. Hence $N(Q_S) = \overline{span\{e_{i+1}, e_{i+2}, \cdots\}}$. It is easy to show that $N(Q_S) \in LatS$. Moreover, we know that the weight sequence of $Q_S^{\frac{1}{2}} S Q_S^{+\frac{1}{2}}$ is $\{(a_2^2-a_1^2)^{\frac{1}{2}}, a_2(\frac{a_3^2-a_2^2}{a_2^2-a_1^2})^{\frac{1}{2}}, \cdots, a_{i-1}(\frac{a_i^2-a_{i-1}^2}{a_{i-1}^2-a_{i-2}^2})^{\frac{1}{2}}, 0, 0, \cdots\}$. So $Q_S^{\frac{1}{2}} S Q_S^{+\frac{1}{2}} \in B(H)$. By the theorem of Section 1, $S$ is a near subnormal operator. Hence Theorem 1 is proved.

### 3. AN APPLICATION

In Theorem 1(2), if $i = 1, 2,$ then it is easy to show that $S$ is subnormal as well as near subnormal by use of the criteria for subnormal operators of [1] or [8].

In the case of $i \geq 3,$ the following theorem holds.

**Theorem 2.** Suppose $S$ is a hyponormal operator with $Se_n = a_n e_{n+1}$ for all $n \geq 1$. If there exists some integer $i \geq 3$ such that $0 < a_1 < a_2 < \cdots < a_{i-1} < a_i$ and $a_n = a_i$ for all $n \geq i+1,$ then $S$ is near subnormal but not subnormal.

Proof. By Theorem 6 of [3] or the theorem of [4], it is not difficult to show that $S$ is not subnormal. Hence Theorem 2 follows directly from Theorem 1.

Now, by Theorem 2, we can construct many operators which are near subnormal but not subnormal. Similarly, from Theorem 1, we can also construct many operators which are hyponormal but not near subnormal. All of these operators can serve as further answers to the Hilbert space problem 160.

Two examples are given at the end of the paper.

**Example 1.** Suppose $a, b$ and $c$ are any numbers such that $0 < a < b < c$. Let $Se_n = a_n e_{n+1}$ for all $n \geq 1$, where $a_1 = a, a_2 = b$ and $a_n = c$ for $n \geq 3$. Then $S$ is near subnormal but not subnormal.

**Example 2.** Let $Se_n = a_n e_{n+1}$ for all $n \geq 1$, where $a_1 = a_2 = k$ ($0 < k < \frac{3}{4}$) and $a_n = 1 - \frac{1}{n+1}$ for $n \geq 3$. Then $S$ is hyponormal but not near subnormal.